\documentclass[]{amsart}

\newtheorem{theorem}{Theorem}[section]
\newtheorem{lemma}[theorem]{Lemma}

\newtheorem{remark}[theorem]{Remark}
\newtheorem{corollary}[theorem]{Corollary}

\newtheorem{proposition}[theorem]{Proposition}

\newtheorem{definition}[theorem]{Definition}

\usepackage{graphicx,graphics}
\usepackage[all]{xy}
\newcounter{rmnum}

\numberwithin{equation}{section}
\begin{document}

  \baselineskip=17pt

\title{Equivariant extension properties of coset spaces of locally compact groups and approximate slices}

\author{Sergey A. Antonyan\\[3pt]}

\address{Departamento de  Matem\'aticas,
Facultad de Ciencias, Universidad Nacional Aut\'onoma de M\'exico,
 04510 M\'exico Distrito Federal, Mexico.}
\email{antonyan@unam.mx}

\begin{abstract} We prove that for a compact subgroup $H$ of a  locally compact Hausdorff group $G$, the following properties are mutually  equivalent: (1) $G/H$ is a manifold, (2)~$G/H$ is  finite-dimensional and locally connected, (3)~$G/H$ is  locally contractible, (4)~$G/H$ is  an ANE for paracompact spaces, (5)~$G/H$ is  a metrizable $G$-ANE for paracompact proper $G$-spaces having a paracompact orbit space.  A  new version of the Approximate slice theorem is also proven  in  the light of these results.

\end{abstract}


\maketitle \markboth{SERGEY ANTONYAN}{EQUIVARIANT EXTENSION PROPERTIES OF COSET SPACES}

\section{Introduction}

In the fundamental work of R.~Palais \cite{pal:60}, it  was established  that for any compact Lie group $G$ and its compact  subgroup $H$ the coset space $G/H$  has the following equivariant extension property: for every normal $G$-space $X$ and a closed invariant subset $A\subset X$, every  $G$-map $f:A\to G/H$ extends to a $G$-map $f' :U\to G/H$ defined on an invariant neighborhood $U$ of $A$. In this case one writes $G/H\in G$-ANE.
This property of $G$ is equivalent to the, so-called, Exact slice theorem:  every orbit in a completely regular $G$-space $X$ is a neighborhood $G$-equivariant retract of  $X$ (see \cite{pal:60}).

 In general, when $G$ is a compact non-Lie group, the Exact slice theorem  is no longer true  (see \cite{ant:94}). At the same time,  among the coset spaces  $G/H$ of a compact non-Lie group $G$, still there are many which  posses the property $G/H\in G$-ANE.
This observation led the author in \cite{ant:94} to the, so-called, Approximate  slice theorem, which is valid  for every compact group of transformations. It claims that given a point $x$ and its neighborhood $O$  in a $G$-space, there exists a  $G$-invariant neighborhood $U$ of $x$ that  admits a $G$-map $f:U\to G/H$ to a coset space $G/H$  with a compact subgroup $H\subset G$ such that  $G/H\in G$-ANE and $x\in f^{-1}(eH)\subset O$.

On the other hand, in 1961 the Exact slice theorem was extended by R.~Palais~\cite{pal:61} to the case of proper actions of non-compact Lie groups. Validity of  the property $G/H\in G$-ANE for every compact subgroup $H$ of a noncompact Lie group $G$ follows from Palais' Exact slice theorem; this was  proved later in  E.~Elfving~ \cite[p.\,23-24]{elf:96}.

It is one of the purposes of this paper to prove (see Proposition~\ref{P:21} and Theorem~\ref{T:0}) that  if $G$ is a locally compact group and $H$ a compact subgroup of $G$ then the following properties are mutually  equivalent: (1) $G/H$ is a manifold, (2)~$G/H$ is  finite-dimensional and locally connected, (3)~$G/H$ is  locally contractible, (4)~$G/H$ is  an ANE for paracompact spaces, (5)~$G/H$ is  a metrizable $G$-ANE for paracompact proper $G$-spaces having a paracompact orbit space.

 One should clarify that in the case of noncompact group actions  the property $G/H\in G$-ANE reads as follows:
 for every paracompact proper $G$-space $X$ having a paracompact orbit space, every  $G$-map $f:A\to G/H$ from  a closed invariant subset $A\subset X$ extends to a $G$-map $f' :U\to G/H$ over  an invariant neighborhood $U$ of $A$.

The equivalence of  the above properties (3) and (5) is claimed also in \cite[Proposition~3.2]{ant:jap}
However,  the proof given in \cite{ant:jap} is valid only for locally compact almost connected groups.
   Proposition~\ref{P:21} and Theorem~\ref{T:0}, in particular,   fill up this gap also.

 As in the compact case \cite{ant:94}, the equivariant extension properties of coset spaces are  conjugated with approximate slices. In section~5 we prove a new version of the Approximate slice theorem valid for proper actions of arbitrary locally compact groups.
Section~3 plays  an auxiliary role. Here we establish   a  special equivariant embedding of a coset space $G/H$ into a  $G$-AE($\mathcal P$)-space (see Proposition~\ref{P:extra}), which is further used  in the proof of  Theorem~\ref{T:0}.

\section {Preliminaries}

Throughout the paper the letter $G$ will denote a locally compact Hausdorff group unless otherwise  stated; by $e$ \ we shall denote the unity of $G$.

All topological spaces and topological groups are assumed to be Tychonoff (= completely regular and Hausdorff).
The basic ideas and facts of the theory of $G$-spaces or topological
transformation groups can be found in G.~Bredon \cite{bre:72} and in  R.~Palais \cite{pal:60}.
Our basic reference on  proper group actions is Palais' article   \cite{pal:61}.  Other good sources are   \cite{abe:78},  \cite{koz:65}.
For equivariant theory of retracts the reader can see, for instance, \cite{ant:87}, \cite{ant:88} and \cite{ant:99}.

For the  convenience of the reader, we recall however  some more  special definitions and facts below.

\medskip

By a $G$-space we mean a triple $(G, X, \alpha)$,   where  $X$ is a topological space, and $\alpha:G\times X\to X$ is a continuous action of the group $G$ on $X$.
 If $Y$ is another  $G$-space, a continuous map $f:X\to Y$ is called a $G$-map or an equivariant map if $f(gx)=gf(x)$ for every $x\in X$ and $g\in G$.
If $G$ acts trivially on $Y$ then we will use the term
\lq\lq invariant map\rq\rq \ instead of \lq\lq equivariant map\rq\rq.

By a normed linear $G$-space (resp., a  Banach $G$-space) we shall mean  a  $G$-space $L$, where $L$
is a normed linear  space (resp., a Banach space) on which $G$ acts by means of {\it linear
isometries},  i.e., $g(\lambda x+\mu y)=\lambda (gx)+\mu (gy)$ and $\Vert gx\Vert=\Vert x\Vert$ for
all $g\in G$, \ $x, y\in L$ and $\lambda$, $\mu\in \Bbb R$.

If $X$ is a $G$-space then for a subset $S\subset X$, \ $G(S)$
denotes the saturation of $S$; i.e., $G(S)$= $\{gs \mid  g\in G,\ s\in S\}$. In particular, $G(x)$  denotes the $G$-orbit $\{gx\in X \mid g\in G \}$ of $x$.  If $G(S)$=$S$ then $S$ is said to be an invariant or $G$-invariant set. The orbit space is denoted by $X/G$.

For a closed subgroup $H \subset G$, by $G/H$ we will denote the $G$-space
of cosets $\{gH \mid  g\in G\}$ under the action induced by left translations.

If $X$ is a $G$-space and $H$  a closed normal subgroup of $G$ then the $H$-orbit space $X/H$  will always be regarded as a  $G/H$-space endowed with the following action of the group $G/H$:
$(gH)*H(x)=H(gx), \ \ \text{where} \ \ gH\in G/H, \ H(x)\in X/H$.

For any $x\in X$, the subgroup   $G_x =\{g\in G \mid  gx=x\}$ is called  the stabilizer (or stationary subgroup) at $x$.

A compatible metric $\rho$ on a $G$-space $X$  is called invariant or $G$-invariant if $\rho(gx, gy)=\rho(x, y)$ for all $g\in G$ and $x, y\in X$.

A locally compact  group $G$ is called {\it almost connected} whenever its quotient space of connected components is compact.

\medskip

 Let  $X$ be  a $G$-space. Two subsets  $U$ and $V$ in  $X$  are called  {\it thin} relative to each other   \cite[Definition 1.1.1]{pal:61}  if the set
$$\langle U,V\rangle=\{g\in G \mid gU\cap V\ne \emptyset\},$$
   called {\it the  transporter} from $U$ to $V$,  has  a compact closure in $G$.

   A subset $U$ of a $G$-space $X$ is called {\it small} \ if  every point in $X$ has a neighborhood thin relative to $U$. A $G$-space $X$
is called  {\it  proper} (in the sense of Palais)  if   every point in  $X$ has a small neighborhood.

Each orbit in a proper $G$-space is closed, and each stabilizer is compact \cite[Proposition~ 1.1.4]{pal:61}.
 Furthermore, if $X$ is a compact  proper $G$-space, then $G$ has to be compact as well.

Important examples of  proper $G$-spaces are the coset spaces $G/H$ with $H$ a compact subgroup of a locally compact group $G$. Other interesting examples the reader can find in  \cite{abe:78}, \cite{ant:98}, \cite{ant:10}, \cite{ill:95} and \cite{koz:65}.
\medskip

In the present  paper we are especially interested in the class $G$-$\mathcal P$ of all paracompact proper $G$-spaces $X$ that  have paracompact orbit space $X/G$. It is a long time standing  open problem whether the orbit space of any paracompact proper $G$-space is paracompact (see \cite{haj:71} and \cite{abe:78}).

 \medskip

A $G$-space  $Y$ is called an equivariant neighborhood extensor  for a given $G$-space $X$   (notation: $Y\in G$-ANE($X$))  if  for any closed invariant subset $A\subset X$  and  any $G$-map $f:A\to Y$,  there exist an invariant neighborhood $U$ of $A$ in $X$ and  a $G$-map $\psi\colon U\to Y$ such that $\psi|_A= f$. If in addition one  always  can take $U=X$, then we say that $Y$ is an equivariant extensor  for $X$ (notation: $Y\in G$-AE($X$)). The map $\psi$ is called a $G$-extension of $f$.

A $G$-space  $Y$ is called an equivariant absolute neighborhood extensor for the class $G$-$\mathcal P$  (notation: $Y\in G$-ANE($\mathcal P$))  if $Y\in G$-ANE($X$) for any $G$-space
$X\in G$-$\mathcal P$.
Similarly, if  $Y\in G$-AE($X$)  for any  $X\in G$-$\mathcal P$, then $Y$ is called an equivariant absolute extensor for the class $G$-$\mathcal P$  (notation: $Y\in G$-AE($\mathcal P$)).

\begin{theorem}[Abels~\cite{abe:78}]\label{T:Ab}
Let  $G$ be a locally compact  group and $T$ a $G$-space such that $T\in K$-{\rm AE}$(\mathcal P)$ for every compact subgroup $K\subset G$. Then $T\in G$-{\rm AE}$(\mathcal P)$.
 \end{theorem}

\begin{remark} In \cite[Theorem~4.4]{abe:78} the result is stated only for  $G$-{\rm AE}$(\mathcal M)$ while the proof is valid also for  $G$-{\rm AE}$(\mathcal P)$, where  $G$-$\mathcal M$ stands for the class of all proper $G$-spaces that are  metrizable by a $G$-invariant metric.
\end{remark}

\begin{corollary} \label{C:B} Let $G$ be a locally compact group and  $B$  a Banach $G$-space. Then $B\in G$-{\rm AE}$(\mathcal P)$.
\end{corollary}
\begin{proof} By a result of E.~Michael~\cite{mich:53}, $B$ is an  AE$(\mathcal P)$. Then, it follows from  \cite[Main Theorem] {ant:80} that $B$ is a $K$-AE$(\mathcal P)$ for every compact subgroup $K$ of $G$. It remains to apply Theorem~\ref{T:Ab}.
\end{proof}

\medskip

Let us recall  the well-known definition of a slice \cite[p.~305]{pal:61}:

\begin{definition}\label{D:21} Let $X$ be a $G$-space and   $H$  be a closed  subgroup of $G$. An $H$-invariant  subset $S\subset X$ is called an  $H$-slice in $X$ if $G(S)$ is open in $X$ and there exists  a $G$-equivariant map $f:G(S)\to G/H$ such that $S$=$f^{-1}(eH)$.  The saturation $G(S)$ is called  a {\it tubular} set.
  If  $G(S)=X$ then we say that $S$ is {\it a global} $H$-slice of $X$.
\end{definition}

\medskip
The following  result of R. Palais \cite[Proposition 2.3.1]{pal:61} plays  a central role in the theory of topological transformation groups:

\begin{theorem}[Exact slice theorem]\label{T:Sl} Let $G$ be a Lie group, $X$  a proper $G$-space and $a\in X$. Then there exists a $G_a$-slice $S\subset X$ such that $a\in S$.
\end{theorem}

 \smallskip

\section{Equivariant embeddings into a  $G$-AE($\mathcal P$)-space}

Recall that the letter $G$ always denotes a   locally compact Hausdorff group.

The main result of this section is Proposition~\ref{P:extra} which provides a  special equivariant embedding of a coset space $G/H$ into a  $G$-AE($\mathcal P$)-space; it is used in the proof of  Theorem~\ref{T:0} below.

We begin with the following lemma proved in \cite[Lemma~2.3]{NSL}:

\begin{lemma} \label{lemma:2}
 Let  $H$  a compact subgroup of $G$ and $X$  a metrizable proper $G$-space
 admitting a  global $H$-slice $S$. Then there is a compatible $G$-invariant metric $d$ on $X$
such that each open unit ball $O_d(x, 1)$    is a small set.
\end{lemma}

\begin{lemma} \label{lemma:3}
 Let   $H$ be   a compact subgroup of $G$.  Then a  subset $S\subset G/H$ is small if and only if the closure $\overline{S}$ is  compact.
\end{lemma}
\begin{proof} Assume that  $S$ is a small subset of $G/H$. Then  there is a neighborhood $U$ of the point $eH\in G/H$ such that the transporter $$\langle S, U\rangle=\{g\in G\mid gS\cap U\ne\emptyset\}$$
has a compact closure. Due to local compactness of $G/H$ one can assume that the closure $\overline{U}$ is compact.

Next, for every $sH\in S$ we have $s^{-1}sH=eH\in U$, i.e., $s^{-1}\in \langle S, U\rangle$, or equivalently, $s\in \langle U, S\rangle$. Hence $sH\in \langle U, S\rangle (U)$ showing that $S\subset  \langle U, S\rangle (U)\subset  \overline{\langle U, S\rangle} (\overline{U})$. But, due  to compactness of the closures   $\overline{\langle U, S\rangle}$ and $\overline{U}$, the set  $\overline{\langle U, S\rangle} (\overline{U})$ is compact. This yields that the closure  $\overline{S}$ is  compact, as required.

The converse is immediate  from the fact that $G/H$ is a proper $G$-space and every compact subset of a proper $G$-space is a small set (see \cite[\S1.2]{pal:61}).

\end{proof}

\begin{corollary} \label{C:1}
 Let  $H$ be a compact subgroup of $G$ such that the quotient $G/H$ is metrizable.  Then there is a compatible $G$-invariant metric $\rho$ on $G/H$ such that each closed unit ball  $B_\rho(x, 1)$, $x\in G/H$,    is compact.
\end{corollary}
\begin{proof} By lemmas\,\ref{lemma:2} and \ref{lemma:3}, there exists a compatible $G$-invariant metric $d$ in $G/H$ such that  each open unit ball $O_d(x, 1)$    has a compact closure. Then the metric $\rho$ defined by  $\rho(x, y)=2d(x, y)$, $x, y\in G/H$, has  the desired property because $B_\rho(x, 1)\subset \overline{O_d(x, 1)}$.
\end{proof}

\medskip

 We recall that a continuous function $f:X \to \mathbb{R}$ defined on a
$G$-space $X$ is called $G$-uniform if for each $\epsilon>0$, there is an identity  neighborhood $U$ in  $G$ such that $|f(gx)-f(x)|< \epsilon$ for all $x \in X$ and $g \in U$.

 For a proper $G$-space $X$ we denote by $\mathcal P(X)$  the linear space of  all bounded  $G$-uniform
functions  $f:X\to \Bbb R$  whose support ${\rm supp}\,
f=\{x\in X \ | \ f(x)\ne 0\}$ is a small subset of $X$.
 We endow $\mathcal P(X)$ with
the sup-norm and the following $G$-action:
$$ (g, f)\mapsto gf, \quad (gf)(x)=f(g^{-1}x), \ \ x\in X.$$

It is easy to see that  $\mathcal P(X)$ is a normed linear $G$-space.  It was  proved in \cite[Proposition\,3.1]{NSL} that the complement  $\mathcal P_0(X)=\mathcal P(X)\setminus\{0\}$ is a proper $G$-space. Moreover, it follows immediately from \cite[propositiones\, 3.4 and 3.5]{ant:FM09} that   the complement  $\widetilde{\mathcal P_0}(X)=\widetilde{\mathcal P}(X)\setminus\{0\}$ is also a  proper $G$-space, where $\widetilde{\mathcal P}(X)$ denotes the completion of $\mathcal P(X)$.

The following  result follows immediately from \cite[Lemma\,3.3]{NSL} and  \cite[proof of Proposition~3.4]{NSL}:

\begin{proposition} \label{prop:4}
Let  $(X,\rho)$ be  a metric proper $G$-space
with an invariant metric $\rho$ such that each closed unit ball in $X$ is
a small set. Then $X$ admits a  $G$-embedding $i:X
\hookrightarrow \mathcal P_0(X)$ such that:
\begin{enumerate}
\item $\Vert i(x)-i(y)\Vert\le \rho(x, y)$ for all $x, y\in X$,
\item $\rho(x, y)=\Vert i(x)-i(y)\Vert$ whenever $\rho(x, y)\le 1$,
\item  $\Vert i(x)-i(y)\Vert \ge1$ whenever $\rho(x, y)>1$,
\item the image  $i(X)$ is closed in its  convex hull.
\end{enumerate}
\end{proposition}

We aim at applying this result to the case  $X =G/H$,  where $H$ is a compact subgroup of $G$ such that $G/H$ is metrizable.

As it follows from Lemma~\ref{lemma:3}, in this specific case,   $\mathcal P(G/H)$ is just the space of all continuous functions $G/H\to  \Bbb R$ having a precompact support. Respectively, $\widetilde{\mathcal P}(G/H)$ is the Banach space of  all continuous functions $G/H\to  \Bbb R$ vanishing at infinity.

\begin{proposition} \label{prop:5} Let  $H$ be  a compact subgroup of $G$ such that the quotient space  $G/H$ is metrizable.
Choose, by Corollary~\ref{C:1}, a compactible $G$-invariant metric $\rho$  on $G/H$ such that each closed unit ball  $B_\rho(x, 1)$, $x\in G/H$,    is compact.
Then $G/H$ admits a $G$-embedding $i:G/H \hookrightarrow \widetilde{\mathcal P_0}(G/H)$ such that:
\begin{enumerate}
\item [(a)] $\Vert i(x)-i(y)\Vert\le \rho(x, y)$ for all $x, y\in G/H$,
\item [(b)] $\rho(x, y)=\Vert i(x)-i(y)\Vert$ whenever $\rho(x, y)\le 1$,
\item [(c)] $\Vert i(x)-i(y)\Vert \ge1$ whenever $\rho(x, y)>1$,
\item [(d)] the image  $i(G/H)$ is closed in $\widetilde{\mathcal P}(G/H)$.
\end{enumerate}
\end{proposition}
\begin{proof} Since $(G/H, \rho)$ satisfies the hypothesis of  Proposition\,\ref{prop:4},  there exists   a topological $G$-embedding $j:G/H
\hookrightarrow \mathcal P_0(G/H)$ satisfying all the four properties  in  Proposition\,\ref{prop:4}. Composing this $G$-embedding with the isometric $G$-embedding   $\mathcal P_0(G/H)\hookrightarrow \widetilde{\mathcal P_0}(G/H)$ we get a  $G$-embedding  $i:G/H
\hookrightarrow \widetilde{\mathcal P_0}(G/H)$ that also satisfies the four properties  in  Proposition\,\ref{prop:4}.
Hence,  the above properties  (a), (b) and (c) are fulfilled.  So, only the last property (d) needs to be verified.

 Let $(x_n)_{n\in \Bbb N}$ be a sequence in $G/H$ such that $\big(i(x_n)\big)_{n\in \Bbb N}$ converges to a point $f\in \widetilde{\mathcal P}(G/H)$. One should to check that $f\in i(X)$.  Since   $\big(i(x_n)\big)_{n\in \Bbb N}$ is a Cauchy sequence, it  follows from the above property (b)  that $(x_n)_{n\in \Bbb N}$ is also Cauchy. Since each closed unit ball  $B_\rho(x, 1)$    is a  compact subset of $G/H$, it then follows from \cite[Ch.\,XIV, Th.\,2.3]{dug:66} that  $\rho$ is a complete metric. Then  $(x_n)_{n\in \Bbb N}$ converges to a limit, say $y\in G/H$. By continuity of $i$, this implies that  $i(x_n)\rightsquigarrow i(y)$, and hence, $f=i(y)\in i(G/H)$, as required.
\end{proof}

\begin{lemma}\label{L:0}  Let  $H$ be a compact subgroup of $G$ and $G_0$ the connected component of  $G$.  If $G/H$ is locally connected then  the subgroup $G_0H\subset G$ is open and almost connected.

\end{lemma}
\begin{proof} Since the natural map
$$G/H\to G/G_0H, \qquad gH\mapsto gG_0H$$
is open and the  local connectedness is invariant under   open maps, we infer that $G/G_0H$ is locally connected. On the other hand
$${G/G_0H}\cong\frac{G/G_0}{(G_0H)/G_0}.$$
Consequently, $G/G_0H$, being the quotient space of the totally disconnected group $G/G_0$ is itself totally disconnected. Hence, $G/G_0H$ should be discrete, implying that $G_0H$ is an open subgroup of $G$.

To prove  that $G_0H$ is almost connected it suffices to observe that the  quotient group $G_0H/G_0$  is just the image of the compact group  $H$ under the natural homomorphism $G\to G/G_0$, and hence, is compact.
\end{proof}

\begin{proposition}\label{P:extra}
Let $H$  be a compact  subgroup of $G$ such that $G/H$ is metrizable and locally connected. Then there exists a closed  $G$-embedding $i:G/H \hookrightarrow \widetilde{\mathcal P_0}(G/H)$ such that
\begin{equation}\label{different}
   \left\{ \begin{array}{rcl}
\Vert i(x)-i(y)\Vert> 1/2  \ \ \text{whenever} \ x\ \text{and} \ y \ \text{ belong}\\ \text{ to different connected components of} \ G/H.
\end{array}\right.
\end{equation}

 \end{proposition}

\begin{proof} Since  $G/H$ is metrizable,  by virtue of Corollary\,\ref{C:1},  there exists a  compatible $G$-invariant metric $d$ on $G/H$  such that each closed unit ball  $B_d(x, 1)$, $x\in G/H$,    is compact.

 To shorten our notation,  set $S=G_0H/H$, where $G_0$ stands for the identity component of $G$.  Since $S$ is the image of $G_0$ under the quotient map $G\to G/H$ we infer that $S$ is connected.
 On the other hand, it follows from Lemma~\ref{L:0} that  $S$ is an open and closed subset of $G/H$, so $S$ should be a  connected component of $G/H$.
   Further, it is easy to check that $S$ is a global $G_0H$-slice for $G/H$. Thus,
   $$G/H= G(S)= \bigsqcup\limits_{g\in G} gS,$$
    the disjoint union of  closed and open connected components  $gS$, one $g$ out of every coset in $G/G_0H$.

   Next we define a new metric $\rho$ on $G/H$ as follows:
      $$
     \left\{ \begin{array}{rcl}
           \text{if two points} \  x \  \text{and} \ y \  \text{ of} \ G/H \ \text{belong to the same connected component,}\\
           \text{then we put}  \             \rho(x, y)=d(x, y); \ \text{otherwise we set} \ \rho(x, y)=d(x, y)+1/2.
           \end{array}\right.
$$

 Clearly, $\rho$ is a compatible metric for $G/H$. Since    $B_\rho(x, 1)\subset B_d(x, 1)$ for every $x\in X$,  we see that each closed unit ball  $B_\rho(x, 1)$ is compact.

 To see the $G$-invariance of $\rho$ assume that $x, y\in X$ and $h\in G$. If the points $x$ and $y$  are  in the same connected component $gS$ then the points  $hx$ and $hy$  belong to the same connected component $hgS$. But then $\rho(hx, hy)= d(hx, hy)=d(x, y)=\rho(x, y)$, as required.

 By the same argument, if $x$ and $y$ are in two different  connected components, then  $hx$ and $hy$ also belong to different  connected components. In this case  $\rho(hx, hy)= d(hx, hy)+1/2 =d(x, y)+1/2=\rho(x, y)$, as required. Thus $\rho$ is $G$-invariant.

Now, let $i:G/H \hookrightarrow \widetilde{\mathcal P_0}(G/H)$ be the  closed  $G$-embedding from Proposition\,\ref{prop:5}.
It then follows from the properties (a) and (b) of Proposition\,\ref{prop:5} that (\ref{different}) is satisfied, as required.
\end{proof}

\smallskip

\section{Large subgroups}

Recall that the letter $G$ always denotes a   locally compact  Hausdorff group unless otherwise stated.

\begin{definition}\label{D:large} A compact subgroup $H$ of  $G$ is called {\it large} if the quotient space $G/H$ is locally connected and finite-dimensional.
\end{definition}

The notion of a large subgroup of a compact group first was singled out  in 1991  by the author  \cite{tiraspol} in form of two other its characteristic properties, namely:  \lq\lq$G/H$ is a manifold\rq\rq \  and \lq\lq$G/H$ is a $G$-ANR\rq\rq.  \
More systematically this notion was studied  later in \cite{ant:94} (for compact groups) and in \cite{ant:jap} (for almost connected groups). In this section we shall investigate the remaining general case of an arbitrary locally compact group.
Large subgroups play a central role also in section 5.

\smallskip

The following result is immediate from Lemma~\ref{L:0}:

\begin{corollary} \label{C:0} Let  $H$ be a large  subgroup of $G$ and $G_0$ the connected component of  $G$.  Then   the subgroup $G_0H\subset G$ is open and almost connected.
\end{corollary}

\begin{proposition}\label{P:0} Let $H$ and $K$ be  compact subgroups of $G$ such that $H\subset K$. If  $H$ is a large subgroup then so is $K$.
\end{proposition}

\begin{proof}
Being $H$ a large subgroup of $G$, the quotient $G/H$ is finite-dimensional and locally connected. Since the map $G/H\to G/K$, $gH\mapsto gK$, is continuous and open, we infer that $G/K$ is locally connected. Its finite-dimensionality follows from  the  one  of $G/H$ and the following  equality  (see \cite[Theorem~10]{skl:64}):
\begin{equation}\label{skl}
\dim G/H= \dim G/K +\dim K/H.
\end{equation}
\end{proof}

\begin{proposition}\label{P:-1}  Let $H$ and  $K$ be two compact subgroups of $G$ such that  $K$ is a large subgroup of $G$ while $H$ is a large subgroup of $K$. Then $H$ is a large subgroup of $G$.
\end{proposition}
 \begin{proof}   Being $K$ a large subgroup of $G$, the quotient $G/K$ is finite-dimensional and locally connected. Then the natural map
 $G/H\to G/K$ is a locally trivial fibration with the fibers homeomorphic to $K/H$ (see \cite[Theorem~$13^\prime$]{skl:64}). But $K/H$ is also locally connected (and finite-dimensional) since $H$ is a large subgroup of $K$. This yields that $G/H$ is locally connected. Finite-dimensionality of $G/H$ follows from  the  one  of $G/K$ and $K/H$, and the above  equality  (\ref{skl}).

Thus, $G/H$ is locally connected and finite-dimensional, and hence, $H$ is a large subgroup of $G$, as required.
\end{proof}

The following lemma is a well-known result in the theory of Lie groups
(see e.g., \cite[Ch.~II, Theorem~4.2]{hel:01}):

\begin{lemma}\label{L:Lie}
Let $\Gamma$ be a  Lie group and $\Delta$ a closed subgroup of \ $\Gamma$.
Then the quotient space $\Gamma/\Delta$ has a unique smooth structure
with the property that the natural  $\Gamma$-action on $\Gamma/\Delta$
induced by left translations is smooth.
\end{lemma}

The following theorem of Montgomery and Zippin (see \cite[\S 6.2]{mozi:55}) plays a key role in what follows:

\begin{theorem}[Montgomery-Zippin]\label{T:MZ}
Let  $G$ be an almost connected group that acts effectively and transitively on a locally compact, locally connected, finite-dimensional space. Then $G$ is a Lie group.
 \end{theorem}

\begin{lemma}\label{L:sigma}
Every almost connected group is $\sigma$-compact.
\end{lemma}
\begin{proof} Let $G$ be an almost connected group. By a well-known Malcev-Iwasawa theorem (see \cite[Ch.\,H, Theorem~32.5] {stroppel}),  $G$ has a maximal compact subgroup $K$,   i.e., every compact subgroup of $G$ is   conjugate to a subgroup of $K$.  According to \cite[Theorem~8.5]{ant:FM09}, $G=KG_0$ where  $G_0$ denotes the connected component of $G$. On the other hand, it is a well-known easy fact that every connected locally compact group is $\sigma$-compact, so $G_0=\bigcup\limits_{n=1}^\infty C_n$, where $C_1,$ $C_2, \dots$, are compact subsets of $G_0$. Then, clearly, $G=KG_0=\bigcup\limits_{n=1}^\infty KC_n$, and since $KC_n$ is compact for every $n=1, 2, \dots$, we conclude that $G$ is $\sigma$-compact.
\end{proof}

The following version of the above Montgomery-Zippin theorem  plays an important  role in the proof of Proposition~\ref{P:21} below:

\begin{theorem}\label{T:szente}
Let  $G$ be an almost connected group that acts effectively and transitively on a locally compact, locally contractible space. Then $G$ is a Lie group.
 \end{theorem}
\begin{proof}
This theorem  is essentially  proved   in J. Szente~\cite{szente}: indeed, since by Lemma~\ref{L:sigma},  every almost connected group is $\sigma$-compact, the assertion  follows from \cite[Theorem~4]{szente}.

\end{proof}

\begin{proposition}\label{P:22}
Let $H$  be a compact   subgroup of an almost connected group $G$. Then the following assertions are  equivalent:
\begin{enumerate}
\item $H$ is a large subgroup,
\item There exists a compact normal subgroup $N$ of $G$ such that $N\subset H$ and $G/N$ is a Lie group. In particular,  $G/H$ is a  coset space of a Lie group.
\item $G/H$ is locally contractible.
\end{enumerate}
 \end{proposition}

\begin{proof}
Denote by $N$ the kernel of the $G$-action on  $G/H$, i.e.,
$$N=\{g\in G  \mid   gt=t  \ \text{for all} \ t\in  G/H \}.$$
 Evidently,  $N\subset H$ and the group $G/N$ acts effectively and transitively on  $G/H$.
Since $G$ is an almost connected group,  by virtue of  \cite[Proposition~8.1]{ant:FM09}, the quotient group $G/N$ is also almost connected.

 $(1)\Longrightarrow (2)$. Since $G/H$ is locally connected and finite-dimensional,  one can apply  Theorem~\ref{T:MZ} according to which  $G/N$ is a Lie group. Since   $H/N$  is a compact subgroup of   $G/N$, by Lemma~\ref{L:Lie}, the coset space  $\frac{G/N}{H/N}$ \ is a smooth manifold.

 It remains to observe that the  following homeomorphism holds:
 \begin{equation}\label{homeo}
G/H\cong\frac{G/N}{H/N}.
\end{equation}

\noindent
$(2)\Longrightarrow (3)$ is evident.

\smallskip

$(3)\Longrightarrow (1)$. As above, denote by $N$ the kernel of the $G$-action on  $G/H$.

 Evidently,  $N\subset H$ and the group $G/N$ acts effectively and transitively on  $G/H$.
 Besides,  since $G$ is an almost connected group,  by virtue of  \cite[Proposition~8.1]{ant:FM09}, the quotient group $G/N$ is also almost connected.

 Consequently, one can apply  Theorem~\ref{T:szente} according to which  $G/N$ is a Lie group. Since   $H/N$  is a compact subgroup of   $G/N$, by Lemma~\ref{L:Lie}, the coset space  $\frac{G/N}{H/N}$ \ is a smooth manifold.

 Due to the    homeomorphism (\ref{homeo}) we get  that $G/H$ is a manifold, and in particular, it is finite-dimensional and locally connected. Thus, $H$ is a large subgroup of $G$.
\end{proof}

\begin{remark} In \cite[Proposition~3.2]{ant:jap} it is claimed that  the  assertions $(2)$ and $(3)$ in this proposition are equivalent even for arbitrary locally compact group $G$. Unfortunately, the proof given in  \cite{ant:jap} contains a gap; in fact it is valid only for almost connected  $G$. The error   occurred because in \cite{ant:jap} Theorem~\ref{T:szente} was  inaccurately applied to  arbitrary locally compact groups while it is stated only  for almost connected ones.

\end{remark}

However, for arbitrary locally compact groups we have the following characterization of large subgroups:

\begin{proposition}\label{P:21}
Let $H$  be a compact   subgroup of $G$. Then the following conditions are  equivalent:
\begin{enumerate}
\item $H$ is a large subgroup,
\item $G/H$ is a smooth manifold; in this case it is the disjoint union of open submanifolds  which  all are  homeomorphic to the same coset space of a Lie group,
\item $G/H$ is locally contractible.
\end{enumerate}
 \end{proposition}

\begin{proof} $(1)\Longrightarrow (2)$. Since, by Corollary~\ref{C:0}, $G_0H$ is open in $G$ we see that $G$ is the disjoint union of the open cosets $xG_0H$, $x\in G$. Since  the quotient map $p:G\to G/H$ is continuous, open and closed we  infer that $G_0H/H$ is open and closed  in $G/H$,  and $G/H$ is the disjoint union of its open subsets  $xG_0H/H$, $x\in G$. Observe that each $xG_0H/H$ is homeomorphic to the coset space $G_0H/H$.

Hence, by virtue of Proposition~\ref{P:22},  it suffices to show that $H$ is a large subgroup of the almost connected group $G_0H$ (see Corollary~\ref{C:0}).

But this is easy. Indeed, since  $G_0H/H$ is an open subset of the  locally connected space $G/H$, we infer that  $G_0H/H$ is locally connected too. Further,  since $G_0H/H$ is closed in  $G/H$ we infer that $\dim  G_0H/H\le \dim  G/H$, and hence, $G_0H/H$ is finite-dimensional because $G/H$ is so. Thus, $H$ is a large subgroup of $G_0H$, as required.

\noindent
$(2)\Longrightarrow (3)$ is evident.

\smallskip

$(3)\Longrightarrow (1)$. Since local contractibility yields  local connectedness, it follows from Corollary~\ref{C:0} that $G_0H$ is an open almost connected subgroup of $G$. In turn, $G_0H/H$ is an open subset of the locally contractible space $G/H$, and hence, is itself locally contractible.
Then we can apply Proposition~\ref{P:22} to the almost connected group $G_0H$ according to which $H$ is a large subgroup of $G_0H$, and hence, the quotient
 $G_0H/H$ is a manifold.

But $G/H$ is the  disjoint union of its open and closed subsets all homeomorphic to $G_0H/H$. This  yields that $G/H$ is a manifold, and in particular, it is finite-dimensional and locally connected. Thus, $H$ is a large subgroup of $G$.
\end{proof}

The following result is proved in Elfving~\cite[p.~23-24]{elf:96} in a different way:

\begin{proposition}\label{P:20} Let $G$ be a Lie group  and  $H$  a compact subgroup of $G$. Then
 $G/H$ is a $G$-{\rm ANE}$(\mathcal P)$.
 \end{proposition}

\begin{proof} By virtue of Corollary\,\ref{C:1},  there exists a  compatible $G$-invariant metric $\rho$ on $G/H$  such that each closed unit ball  $B_\rho(x, 1)$, $x\in G/H$,    is compact.

Next, by Proposition\,\ref{prop:5}, one can assume that $G/H$ is an invariant closed subset of the proper $G$-space
 $\widetilde{\mathcal P_0}(G/H)$. Now, due to  Exact slice theorem~\ref{T:Sl} (see also \cite[Corollary~1]{pal:61}), $G/H$ is a  $G$-retract of some invariant neighborhood $U$ in  $\widetilde{\mathcal P_0}(G/H)$. Since  by Corollary~\ref{C:B},  $\widetilde{\mathcal P}(G/H)\in G$-AE($\mathcal P$), we conclude that $G/H\in G$-ANE($\mathcal P$), as required.

\end{proof}

\begin{proposition}\label{P:Almost} Let $G$ be an almost connected group and  $H$  a large   subgroup of $G$. Then
 $G/H$ is a metrizable $G$-{\rm ANE}$(\mathcal P)$.
 \end{proposition}

\begin{proof}  Denote by $N$ the kernel of the $G$-action on  $G/H$, i.e.,
$$N=\{g\in G  \mid   gt=t  \ \text{for all} \ t\in  G/H \}.$$
 Then  $N\subset H$ and the group $G/N$ acts effectively and transitively on the finite-dimensional locally connected space $G/H$.  Consequently, according to  Theorem~\ref{T:MZ}, $G/N$ is a Lie group.

 We  have  to return again to the  $G$-equivariant homeomorphism (\ref{homeo}).

Since $G/N$ is a Lie group, it then follows from  Proposition~\ref{P:20} and homeomorphism (\ref{homeo}) that   $G/H$ is a $G/N$-ANE$(\mathcal P)$. Further, $G/H$ is metrizable since it is homeomorphic to the coset space of the metrizable (in fact Lie)  group $G/N$.

 Now, since $N$ acts trivially on $G/H$, it then follows  from \cite[Proposition 3]{ant:87} that $G/H$ is a $G$-ANE$(\mathcal P)$.
\end{proof}

We have developed all the tools necessary to prove our  main result:

\begin{theorem}\label{T:0}
Let $H$  be a compact  subgroup of $G$. Then the following conditions are  equivalent:
\begin{enumerate}
\item $H$ is a large subgroup,
\item $G/H$ is a metrizable $G$-{\rm ANE}$(\mathcal P)$,
\item $G/H$ is an {\rm ANE}$(\mathcal P)$.

\end{enumerate}
 \end{theorem}
\begin{proof} $(1)\Longrightarrow (2)$. By Proposition~\ref{P:21}, $G/H$ is metrizable. Then, by Proposition~\ref{P:extra}, one can assume that $G/H$ is a $G$-invariant closed subset of  $\widetilde{\mathcal P_0}(G/H)$ and
\begin{equation}\label{different1}
   \left\{ \begin{array}{rcl}
\Vert x-y\Vert> 1/2  \ \ \text{whenever} \ x\ \text{and} \ y \ \text{ belong}\\ \text{ to different connected components of} \ G/H.
\end{array}\right.
\end{equation}

Set $S=G_0H/H$ and denote by $W$ the $1/4$-neighborhood of $S$ in $\widetilde{\mathcal P_0}(G/H)$, i.e.,
$$W=\{z\in \widetilde{\mathcal P_0}(G/H)\mid dist(z, S)<1/4\}.$$

\noindent
{\it Claim}. \ $W$ is a $G_0H$-slice in  $\widetilde{\mathcal P_0}(G/H)$.

Indeed, $W$ is $G_0H$-invariant since $S$ is so and the norm of $\widetilde{\mathcal P}(G/H)$ is $G$-invariant. Further, $G(W)$ is open in $\widetilde{\mathcal P_0}(G/H)$ since $W$ is so.

Check that $W$ and $gW$ are disjoint whenever $g\in G\setminus G_0H$. In fact, if $gw\in W$ for some $w\in W$ then
$ \Vert w-s\Vert<1/4$ and $ \Vert gw-s_1\Vert<1/4$ for some $s, s_1\in S$. By the invariance of the norm we have $ \Vert gw-gs\Vert =  \Vert w-s\Vert<1/4$. Hence,
$$ \Vert gs-s_1\Vert\le \Vert gs-gw\Vert  + \Vert gw-s_1\Vert < 1/4+1/4=1/2.$$
Consequently, by (\ref{different1}), $s_1$ and $gs$ must belong to the same connected component of $G/H$. Since $s_1\in S$  we infer that $gs\in S$. Thus, $S\cap gS\ne\emptyset$.
But, since $S$ is a global $G_0H$-slice of  $G/H$ (see the proof of Proposition~\ref{P:extra}), it then follows  that $g\in G_0H$, as required.

Thus, $G(W)$ is the disjoint union of its open subsets $gW$. In particular, each $gW$ is also  closed in  $G(W)$.

So, we have verified that $W$ is a $G_0H$-slice in   $\widetilde{\mathcal P_0}(G/H)$.

Further, since $G_0H$ is an almost connected group (see Corollary~\ref{C:0}), it follows from   Proposition~\ref{P:Almost} that   $S$ is a $G_0H$-ANE$(\mathcal P)$. Now, since  $W\in G_0H$-$\mathcal P$, it then follows  that there exists a $G_0H$-equivariant retraction $r:V\to S$ for  some $G_0H$-invariant  neighborhood $V$  of $S$ in $W$.

Then $r$ induces  a $G$-map $R:G(V)\to G/H$ by the rule: $R(gv)=gr(v)$, where $g\in G$ and $v\in V$ (see \cite[Ch.~I, Proposition~ 4.3]{tdie:87}).

Clearly, $R$ is a $G$-retraction, and hence, $G/H$ being a $G$-neighborhood retract of the $G$-AE($\mathcal P$)-space $\widetilde{\mathcal P}(G/H)$ (see Corollary~\ref{C:B}), is itself a $G$-ANE($\mathcal P$)-space.

\smallskip

$(2)\Longrightarrow (3)$. Suppose that $X$ is a paracompact space,  $A$ a closed subset of $X$, and $f:A\to G/H$  a continuous map. Consider the $G$-space $G\times X$ endowed with the action of $G$ defined by the rule:
$h(g, x)=(hg, x)$ for all $(g, x)\in G\times X$ and $h\in G$. Then the map $F:G\times A\to G/H$, given by $F(g, a)= gf(a)$, is a continuous $G$-map.

Since $G$ is a proper $G$-space,  so is the product $G\times X$.
Since the orbit space of $G\times X$ is evidently homeomorphic to $X$, we conclude that $G\times X\in G$-$\mathcal P$.
  Hence, by the hypothesis, there exist a $G$-neighborhood $U$ of $G\times A$ in $G\times X$ and a $G$-map $\widetilde F: U\to G/H$ that extends $F$. Since $\{e\}\times A\subset U$ we infer that $\{e\}\times V\subset U$ for some neighborhood $V$ of $A$ in $X$. But $U$ is a $G$-invariant set, implying that  $G\times V\subset U$.

Next we define a map $\widetilde f:V\to G/H$ by putting $\widetilde f(v)=\widetilde F(e, v)$. Clearly, $\widetilde f$ \ is a continuous extension of $f$, as required.

\smallskip

$(3)\Longrightarrow (1)$. By Proposition~\ref{P:21}, it suffices to show that  $G/H$ is locally contractible. Assume that $U$ is a compact neighborhood of the point $a=eH$ in $G/H$. Denote by  $A$  the closed subset $U\times\{0\}\cup \{a\}\times I\cup U\times\{1\}$ of the product $U\times I$, where $I=[0, 1]$, and consider the continuous map $f:A\to G/H$ defined by the rule:
$$f (u,0)=u  \ \ \text{and} \ \ f(u, 1)=a \ \  \text{if}\  u\in U\quad \text{and}\quad f(a, t)=a \ \ \text{for all} \ \ t\in I.$$
Since $U\times I$ is paracompact (in fact, compact) and $G/H\in {\rm ANE}(\mathcal P)$, the map $f$  extends  to a continuous map  $F:V\to G/H$ over  an open neighborhood $V$ of $A$ in $U\times I$.

Choose an open neighborhood $W$ of $a$ in $U$ such that $W\times I\subset V$. Then the restriction $F|_{W\times I}  : W\times I\to U$ is a contraction of $W$ in $U$ to the point $a$, as required.
\end{proof}

We conclude this section with the following two corollaries, which  were  used in the proof of \cite[Theorem~1.1]{ant:05}:

\begin{corollary}[\cite{ant:05}, Lemma 2.5]\label{C:SE} Let  $H$ be a  large subgroup of $G$. Assume that $A$ is a closed invariant subset of a proper $G$-space $X\in G$-$\mathcal P$,  and  $S$ is a global $H$-slice of $A$. Then there exists an $H$-slice $\widetilde{S}$ in $X$ such that $\widetilde{S}\cap A=S$.
\end{corollary}
\begin{proof} Let $f:A\to G/H$ be a  $G$-map with $f^{-1}(eH)=S$. By Theorem~\ref{T:0}, $G/H\in G$-ANE($\mathcal P$), and hence, there exists a $G$-extension $F:U\to G/H$ over an invariant neighborhood $U$ of $A$ in $X$. It is easy to see that the preimage $\widetilde{S}=F^{-1}(eH)$ is the desired $H$-slice.
\end{proof}

\begin{corollary}[\cite{ant:05}, Lemma 3.2]\label{L:L}  Let  $H$ be a closed normal subgroup of $G$, and $K$ a large subgroup of $G$.  Then  $(K H)/H$  is  a  large subgroup of $G/H$.
\end{corollary}

\begin{proof} Since $(KH)/H$ is the  image of the compact subgroup $K$ under the continuous homomorphism $G\to G/H$ we see that it is a compact subgroup of $G/H$. Further, observe that the following homeomorphism (even $G$-equivariant) holds:
\begin{equation}\label{Skl}
\frac{G/H}{(K H)/H}\cong G/KH.
\end{equation}
Since $K$ is a large subgroup and $K\subset KH$, it follows from  Proposition~\ref{P:0} that $KH$ is so. Hence,  the coset space $G/KH$ is finite-dimensional and locally connected. It remains to apply  (\ref{Skl}).
\end{proof}

\smallskip

\section {Approximate slices for proper actions of non-Lie groups}
\smallskip

In \cite{abe:78} and  \cite{ant:jap}  approximate  versions of the Exact slice theorem \ref{T:Sl} were established which are  applicable also to proper actions of  non-Lie groups.

In this section we shall prove the following  new version of the Approximate  slice theorem for proper actions of {\it arbitrary locally compact groups} which improves the one in \cite[Theorem~3.6]{ant:jap}:

\begin{theorem}[Approximate slice theorem] \label{T:331}  Assume that  $X$ is  a proper $G$-space,   $x\in X$ and $O$  a neighborhood of $x$. Denote by $\mathcal N(x, O)$ the set of all     large subgroups $H$ of $G$ such that $G_x\subset H$ and $H(x)\subset O$.
Then:
\begin{enumerate}
\item   $\mathcal N(x, O)$ is not empty.
 \item For every  $K \in  \mathcal N(x, O)$,  there exists a  $K$-slice  $S$ with   $x\in S\subset O$.
\end{enumerate}
\end{theorem}

In the proof of this theorem we shall need the following lemma:

\begin{lemma}\label{L:open} Let $X$ be a proper $G$-space, $H$ a compact  subgroup of $G$, and $S$ a global $H$-slice of
$X$. Then the restriction $f:G\times S\to X$  of the action is an open map.
\end{lemma}

\begin{proof} Let $O$ be an open subset of $G$ and $U$ be an open subset of $S$.  It suffices to show that the set $OU=\{gu~|~g\in O, \ u\in U\}$ is open in $X$.

Define $W=\bigcup\limits_{h\in H}(Oh^{-1})\times (hU)$.
We claim that
\begin{equation}\label{open}
X\setminus OU=f\bigl((G\times S)\setminus W\bigr).
\end{equation}

Indeed, since $OU=f(W)$ \ and \ $X=f(G\times S)$, the inclusion  $X\setminus OU\subset f\bigl((G\times S)\setminus W\bigr)$ follows.

Let us establish  the converse inclusion
$f\bigl((G\times S)\setminus W\bigr)\subset X\setminus OU$.

Assume the contrary, that there exists a point  $gs\in f\bigr(G\times S)\setminus W\bigr)$ with  $(g, s)\in (G\times S)\setminus W$ such that  $gs\in OU$. Then $gs=tu$ \ for some \ $(t, u)\in O\times U$. Denote $h=g^{-1}t$. Then one has:
 $$s=g^{-1}tu=hu\qquad\text{and}\qquad (g, \ s)=(tt^{-1}g, \ g^{-1}tu)=(th^{-1}, \ hu)\in (O h^{-1})\times (hU).$$
  Since both $s$ and $u$ belong to $S$,  and  $s=hu$, we conclude that  $h\in H$. Consequently, $(O h^{-1})\times (hU)\subset W$, yielding that    $(g, s)\in W$, a contradiction. Thus, the equality (\ref{open})  is proved.

Being a global $H$-slice,  the set $S$ is a closed small subset of $X$. Consequently, by virtue of  \cite[Proposition~1.4]{abe:78},    the restriction of the action  map  $G\times S\to X$ is  closed. Then, since   $(G\times S)\setminus W$  is  a closed  subset of   $G\times S$, the image  $f\bigl((G\times S)\setminus W\bigr)$  is closed in $X$. Finally, together with the equality (\ref{open}), this implies that  $OU$ is open in $X$, as required.

\end{proof}

\smallskip

\noindent
{\it Proof of Theorem~\ref{T:331}}. We first consider two special cases, namely $G$ totally disconnected and $G$  almost connected, and combine the two to get the general result.

Consider the set $V=\{g\in G \mid gx\in O\}$  which is  an open neighborhood of the compact subgroup $G_x$ in $G$.

\smallskip

{\it Case 1}. Let $G$ be totally disconnected. Then there exists  a compact open subgroup $H$ of $G$ such that $G_x\subset H\subset V$ (see \cite[Ch.\,II, \S\,2.3]{mozi:55}). Therefore $G/H$ is discrete, and hence, $H$ is a large subgroup of $G$. Thus,  $H \in  \mathcal N(x, O)$.

Now assume that $K\in\mathcal N(x, O)$. Then $K$ is a compact open subgroup of $G$ (see, e.g., Corollary~\ref{C:0}).
 %is any large subgroup of $G$ such that $G_x\subset K$ and $K(x)\subset O$, or  equivalently,   $K\subset V$. Then $K$ is an open subgroup of $G$.
  Since  $K(x)\subset O$, there exists a neighborhood $Q$ of $x$ such that $KQ\subset O$. Since $K$ is open, by \cite[Proposition 1.1.6]{pal:61}, there exists a neighborhood $W$ of the point $x$ in $X$ such that $\langle W, W\rangle\subset K$. Then the set $S=K(Q\cap W)$ is a $K$-invariant neighborhood of $x$  with $S\subset O$  and
  $\langle S, S\rangle =K^{-1}\langle Q\cap W, Q\cap W\rangle K=K$. Now, the saturation   $U=G(S)$ is the  disjoint union of open subsets $gS$, one $g$ out of every coset in  $G/K$. So, the map $f:U\to G/K$ with $f(u)=gK$ if $u\in gS$, is a well-defined $G$-map  and  $f^{-1}(eK)=S$. Since  $x\in S\subset O$, we are done.

\smallskip

{\it Case 2}. Let $G$ be almost connected.  By compactness of $G_x$, there exists a  unity neighborhood $V_1$ in   $G$  such that $V_1\cdot G_x\subset V$. By a result of Yamabe (see \cite[Ch.\,IV, \S\,46]{mozi:55} or  \cite[Theorem 8]{gl:60}), $V_1$ contains   a compact normal subgroup $N$ of $G$ such that $G/N$ is a Lie group; in particular, $N$ is a large subgroup of $G$. Setting $H=N\cdot G_x$ we get a compact subgroup $H$ of $G$ such that $G_x\subset H\subset V$. Since $N\subset H$ and $N$ is a large subgroup, it follows from Proposition~\ref{P:0}  that $H$ is also a large subgroup. Thus, $H\in \mathcal N(x, O)$.

\smallskip

Now assume that $K\in \mathcal N(x, O)$.
% is any large subgroup of $G$ such that  $G_x\subset K$ and $K(x)\subset O$, or equivalently, $K\subset V$.
Denote by $M$ the kernel of the $G$-action on  $G/K$, i.e.,
$$M=\{g\in G  \mid   gt=t  \ \text{for all} \ t\in  G/K \}.$$
 Then  $M\subset K$ is a compact normal subgroup of $G$ and the group $G/M$ acts effectively and transitively on the finite-dimensional locally connected space $G/K$.  Consequently, according to  Theorem~\ref{T:MZ}, $G/M$ is a Lie group.

Since $K$ is compact and $K(x)\subset O$, there exists a $K$-invariant neighborhood $Q$ of $x$ such that $Q\subset O$.
  Let $p:X\to X/M$ be  the $M$-orbit map. Then   $X/M$ is a proper $G/M$-space \cite[Proposition 1.3.2]{pal:61}, and it is easy to see that  the $G/M$-stabilizer of the point $p(x)\in X/M$ is just the group $K/M$. Now, by  the  Exact slice theorem \ref{T:Sl}, there exists an invariant neighborhood $\widetilde U$ of $p(x)$ in  $X/M$  and a $G/M$-equivariant map
$$\widetilde f:  \widetilde U\to \frac{G/M}{K/M}$$
 such that $\widetilde f(p(x))=K/M$.

Next we shall consider $X/M$ (and its invariant subsets) as a $G$-space endowed with the action of $G$ defined by the natural homomorphism $G\to G/M$. In particular, $\widetilde U$ is a $G$-space.

 Since the two $G$-spaces  $\frac{G/M}{K/M}$ \ and \  $G/K$ are naturally $G$-homeomorphic, we can consider $\widetilde f$ as a $G$-equivariant map from $\widetilde U$ to $G/K$ with $\widetilde f(p(x))=eK$.

Let $S={\widetilde f}^{-1}(eK)$ and $S_1=S\cap p(Q)$. Then $S$ is a global $K$-slice for $\widetilde U$ and $S_1$ is an open $K$-invariant subset of $S$.

 We claim that the  $G$-saturation $U_1 =G(S_1)$ is a tubular set with $S_1$ as a $K$-slice. Indeed, the openness of  $U_1$  in $\widetilde U$, and hence in $X/M$, follows from  Lemma~\ref{L:open}.

To prove  that $S_1$ is a global $K$-slice of $U_1$ it suffices to show  that  ${f_1}^{-1}(eK) = S_1$, where   $f_1:U_1\to G/K$ is the restriction $\widetilde f|_{U_1}$.

 To this end, choose  $x\in {f_1}^{-1}(eK)$ arbitrary. Since  ${f_1}^{-1}(eK)=   S\cap {U_1}$ then  $x=gs_1$ for some   $s_1\in S_1$ and $g\in G$.
  Hence,   $gs_1\in S\cap gS$, which implies  that $g\in K$. Since $S_1$ is $K$-invariant we  infer  that  $x= gs_1\in  S_1$.  Thus,
  ${f_1}^{-1}(eK) \subset S_1$. The converse inclusion   $S_1\subset {f_1}^{-1}(eK)$ is evident, so we get the desired equality ${f_1}^{-1}(eK) = S_1$.

 Thus,   $S_1$ is a  $K$-slice  lying  in $p(Q)$ and containing the point $p(x)\in X/M$.

Now we set  $U=p^{-1}(U_1)$,   $S=p^{-1}(S_1)$, and let $f:U\to G/K$ be the composition ${f_1} p$. Since $S=f^{-1}(eK)\subset p^{-1}(p(Q))=Q\subset O$ and $x\in S$, we conclude that $S$  is the desired $K$-slice.

\smallskip
{\it Case 3}. Let $G$ be arbitrary. First we show that $\mathcal N (x, O)\ne \emptyset$.

Denote by   $G_{0}$   the identity  component of
   $G$ and  $\widetilde G$=$G/G_0$. Set $\widetilde X=X/G_{0}$ and  let $p: X \to \widetilde X$ be the $G_0$-orbit map. Then $\widetilde  X$ is a proper $\widetilde G$-space \cite[Proposition 1.3.2]{pal:61}, and  the stabilizer ${\widetilde G}_{p(x)}$ of the point $p(x)\in X/M$ in $\widetilde G$ is just the group $(G_0\cdot G_x)/G_0$.

Since  $\widetilde G$ is  totally disconnected,  there exists  a compact open subgroup $M$ of $\widetilde G$ such that
 ${\widetilde G}_{p(x)}\subset M$  (see \cite[Ch.\,II, \S\,2.3]{mozi:55}).

 Denote by $\pi:G\to \widetilde G$ the natural homomorphism and let $L=\pi^{-1}(M)$. Then $L$ is a  closed-open subgroup of $G$. Since, clearly, the quotient group  $G/G_0$ is topologically isomorphic to the compact group $M$ we infer that $L$ is  almost connected.
 Hence, we  can apply the fist part of case~2  to   the  almost connected group $L$, the proper $L$-space $X$ and   the neighborhood  $O\subset X$ of the point $x\in X$. Accordingly, there exists a large subgroup $N$ of $L$ such that $L_x\subset N$ and $N(x)\subset O$.

We claim that $N\in \mathcal N (x, O)$. Indeed, since  $(G_0\cdot G_x)/G_0= {\widetilde G}_{p(x)}\subset M$ we infer that $G_0\cdot G_x\subset L$. In particular, this yields that $G_x=L_x$, and hence, $G_x\subset N$. It remains to check that $N$ is a large subgroup of $G$. In fact, since $N$ is a large subgroup of $L$, due to Proposition~\ref{P:21},  the quotient $L/N$ is locally contractible. But $G/N$ is the disjoint union of its  open subsets  of the form $xL/N$, $x\in G$, each of which is homeomorphic to $L/N$. Consequently, $G/N$ is itself locally contractible, and again by  Proposition~\ref{P:21}, this yields that $N$ is a large subgroup of $G$. Thus, we have proved that $N\in \mathcal N (x, O)$, as required.

  Next, we assume that    $K\in \mathcal N (x, O)$. Since $K$   is a large subgroup of $G$, by Corollary~\ref{C:0},  $H=G_0K$ is an open almost connected subgroup of  $G$. Hence  $\widetilde H=G_0K/G_0$ is a compact open subgroup  of $\widetilde G$.
  The inclusion  $G_x\subset K$ easily implies that ${\widetilde G}_{p(x)}\subset \widetilde H$. Respectively, the inclusion $K(x)\subset O$ yields that $\widetilde H \subset p(O)$.  Then, according  to case~1, there exists a $\widetilde G$-map $f_{1}:U_1\to \widetilde G/\widetilde H$  of an open $\widetilde G$-invariant neighborhood  $U_{1}$ of $p(x)$ in $\widetilde X$  to the discrete $\widetilde G$-space    $\widetilde G/\widetilde H$  with $p(x)\in f_1^{-1}(\widetilde e \widetilde H)\subset p(O)$.

\vspace{0.1cm}

%Let $\pi :  G\to\widetilde G$ be the canonical homomorphism and  $H=\pi^{-1}(\widetilde H)$. Then $H$ is an open, almost connected subgroup of $G$.

The inverse image $W_{1}=f_{1}^{-1}(\widetilde e\widetilde H)$ is an open $\widetilde H$-invariant subset of $\widetilde X$; so the set $W=p ^{-1}(W_{1})$ is an open $H$-invariant subset of $X$ with  $x\in W$, $G_{x}\subset K\subset H$ and $K(x)\subset W\cap O$. Since $H/K$ is open in $G/K$ we infer that $K$  is a large subgroup of $H$ (for instance, by Proposition~\ref{P:21}).

Hence, we can  and do apply case~2 of this  proof to   the  almost connected group $H$, the proper $H$-space $W$,  the neighborhood  $O\cap W$ of the point $x\in W$ and the large subgroup $K$ of $H$. Then there exist  an $H$-neighborhood $U$ of $x$ in $W$ and  an $H$-map $f_{0}:U\to H/K$ with $x\in f_{0}^{-1}(eK)\subset O\cap W$.  Next we want to extend $f_{0}$ to a $G$-map $f:G(U)\to G/K$. Since $H/K\subset G/K$, we simply define $f(gu)=gf_{0}(u)$ for $g\in G,\;u\in U$.

It is easy to check that $f$ is a well-defined $G$-map.  Thus,  the $K$-slice $S=f^{-1}(eK)$ is the desired one.
\qed

\medskip

If  $G$ is a Lie group then, clearly, each compact subgroup of $G$ is large. So, in this case,  Theorem~ \ref{T:331} has the following simpler form:

\begin{corollary} Assume that $G$ is a Lie group,  $X$   a proper $G$-space,   $x\in X$ and $O$  a neighborhood of $x$.
Then for each  compact   subgroups $K$ of $G$ such that $G_x\subset K$ and $K(x)\subset O$, there
 exists a  $K$-slice  $S$ such that    $x\in S\subset O$.
\end{corollary}

We derive   from  Theorem~\ref{T:331} yet another corollary applicable to the, so-called, {\it rich} $G$-spaces.

Recall    that a $G$-space $X$ is called {\it  rich},  if for any point $x\in X$ and any its neighborhood  $U\subset X$, there exists a point $y\in U$ such that the stabilizer $G_{y}$ is a large subgroup of $G$ and $G_{x}\subset G_{y}$ (see \cite{ant:94}, \cite{ant:99}).

\begin{corollary} \label{T:441}  Assume that  $X$ is  a rich  proper $G$-space,   $x\in X$ and $O$  a neighborhood of $x$.
Then  there exist a point $y\in O$ with a large stabilizer  $G_y$ containing  $G_x$, and   a $G_y$-slice $S$ such that  $x\in S\subset O$.

\end{corollary}

\begin{proof}
Choose a neighborhood $O'$ of $x$ such that $G_y(x)$ is contained in $O$ for all $y\in O'$ (see \cite[Lemma~3.9]{ant:jap}). Further, since $X$ is a rich proper $G$-space, we can  choose a point $y\in O'$ such that $G_y$ is a large subgroup of $G$ and  $G_x\subset G_y$. Thus, $G_y\in \mathcal N(x, O)$, the set  defined  in Theorem~\ref{T:331}.

If $G(y)=G(x)$ then the stabilizer $G_x$, being conjugate to $G_y$,  is also a large subgroup, and clearly, $G_x\in \mathcal N(x, O)$.
Next, we apply item (2) of Theorem~\ref{T:331} to  $K=G_x$; the resulting   $K$-slice  $S$  is the desired one.

If $G(y)\ne G(x)$ then we first  choose  (due to  \cite[Proposition~1.2.8]{pal:61}) disjoint invariant neighborhoods $A_x$ and $A_y$ of $G(x)$ and $G(y)$, respectively.
Next, we  apply  twice the first assertion of the statement  (2) of    Theorem~\ref{T:331}: first, to   $x\in O\cap A_x$ and     $K=G_y$, and then  to  $y\in O\cap A_y$ and   $K=G_y$. As a result we get two  $G_y$-slices $S_x\subset O\cap A_x$ and $S_y\subset O\cap A_y$ which contain  the points $x$ and $y$, respectively.  Since $A_x\cap A_y=\emptyset$  the union $S=S_x \cup S_y$ is the desired $G_y$-slice.
\end{proof}

In conclusion we show that there are sufficiently many rich $G$-spaces. Indeed, it was proved in \cite{ant:94} that if $G$ is  a compact group then every metrizable $G$-ANE$(\mathcal M)$ is a rich $G$-space, where $G$-$\mathcal M$ stands for the class of all proper $G$-spaces that are metrizable by a $G$-invariant metric. The same was proved in \cite[Proposition~3.10]{ant:jap} for proper actions of almost connected groups. Below we show that it  is  true also for proper actions of   arbitrary locally compact groups.

\begin{proposition}\label{L:20} Every metrizable proper $G$-{\rm ANE}$(\mathcal M$) is a rich $G$-space.
\end{proposition}
 \begin{proof} Let $X$ be a metrizable proper $G$-{\rm ANE}$(\mathcal M$), $x\in X$ and $O$  a neighborhood of $x$. Then by Theorem~\ref{T:331}, there is a compact large subgroup $K\subset G$, and a $K$-slice $S$ such that $x\in S\subset O$. The tube $G(S)$, being an open subset of $X$,  is a $G$-{\rm ANE}$(\mathcal M$) as well. This yields that   $G(S)$  is  a $K$-{\rm ANE}$(\mathcal M$) (see \cite[Proposition~3.4]{ant:jap}). Hence, according to \cite[Proposition (2)]{ant:94}, $G(S)$ is a rich $K$-space, so there exists a point $y\in G(S)\cap O$ such that $K_x\subset K_y$, and $K_y$ is a large subgroup of  $K$.
 It remains  to show that $G_y$ is a large subgroup of $G$.

 First, it follows from  Proposition~\ref{P:-1} that   $K_y$ is a large subgroup of  $G$. Next, since the point  $y$ belongs to the  $K$-slice $S$ we infer that  $K_y=G_y$, and therefore,  $G_y$ is a  large subgroup of   $G$, as required.
\end{proof}

\subsection* {Acknowledgement}

This research  was supported in part by grants IN102608 from PAPIIT (UNAM) and  79536 from CONACYT (Mexico).

\medskip

\bibliographystyle{amsplain}

\end{document}